\documentclass[letterpaper,onecolumn,draftcls]{ieeeconf}
\usepackage[LGR,T1]{fontenc}
\usepackage[latin9]{inputenc}
\usepackage{amsmath}
\usepackage{amssymb,amsfonts}

\usepackage{mathtools}
\mathtoolsset{showonlyrefs=true}
\usepackage{color}

\usepackage[noadjust]{cite}
\newtheorem{corollary}{Corollary}
\newtheorem{lemma}{Lemma}
\newtheorem{theorem}{Theorem}
\newtheorem{assumption}{Assumption}

\newcommand{\R}{\mathbb{R}}

\IEEEoverridecommandlockouts

\title{\LARGE \bf Totally Asynchronous Distributed Quadratic Programming with \\ Independent Stepsizes and Regularizations}
\author{Matthew Ubl$^{\star}$ and Matthew T. Hale$^{\star}$
\thanks{$^\star$Department of Mechanical and Aerospace
Engineering, University of Florida, Gainesville, FL, 32611, USA. Emails: \texttt{\{m.ubl,matthewhale\}@ufl.edu}}}

\begin{document}
\maketitle

\begin{abstract}
Quadratic programs arise in robotics, communications, smart grids,
and many other applications. As these problems grow in size, finding
solutions becomes much more computationally demanding, and new algorithms
are needed to efficiently solve them. Targeting large-scale problems,
we develop a multi-agent quadratic programming framework in which
each agent updates only a small number of the total decision variables
in a problem. Agents communicate their updated values to each other,
though we do not impose any restrictions on the timing with which
they do so, nor on the delays in these transmissions. Furthermore,
we allow weak parametric coupling among agents, in the sense that
they are free to independently choose their stepsizes, subject to
mild restrictions. We show that these stepsize restrictions depend
upon a problem's condition number. We further provide the means
for agents to independently regularize the problem they solve, thereby
improving condition numbers and, as we will show, convergence properties,
while preserving agents' independence in selecting parameters. Simulation
results are provided to demonstrate the success of this framework
on a practical quadratic program.
\end{abstract}

\section{Introduction}
Convex optimization problems arise in a diverse array
of engineering applications, including signal processing~\cite{luo06},
robotics~\cite{schulman14,mellinger11}, 
communications~\cite{chiang05}, 
machine learning~\cite{shalev12},
and many others~\cite{boyd04}. In all of these areas, 
problems can become very large as the number
of network members (robots, processors, etc.) becomes
large.
Accordingly, there has arisen interest in solving
large-scale optimization problems. 
A common feature of large-scale solvers is that they are parallelized
or distributed among a collection of agents in some way. 
As the number of agents grows,
it can be difficult or impossible to ensure synchrony
among distributed computations and communications, and
there has therefore arisen interest in distributed
asynchronous optimization algorithms. 

One line of research considers asynchronous optimization
algorithms in which agents' communication topologies
vary in time. 
A representative sample of this work
includes~\cite{chen12,zhu12,nedic10,nedic09,ram10,lobel11}, 
and these algorithms
all rely on an underlying averaging-based update law, i.e.,
different agents update the same decision variables and then
repeatedly average their iterates to mitigate disagreements
that stem from asynchrony. 
These approaches (and others in the literature)
require some form of graph connectivity over intervals of a
finite length. 
In this paper, we are interested in cases in which delay bounds are
outside agents' control, e.g., 
due to environmental hazards and adversarial jamming for
a team of mobile autonomous agents. 
In these settings,
verifying graph connectivity can be difficult for single
agents to do, and it may not be possible to even check that
connectivity assumptions are satisfied. 
Furthermore, even if such checking is possible, it will
be difficult to reliably attain connectivity with
unreliable and impaired communications. These
challenges can cause existing approaches to perform poorly
or fail outright. For multi-agent systems with
impaired communications, we are interested in developing an
algorithmic framework that succeeds without requiring
any delay bound assumptions. 

In this paper, we develop a totally asynchronous quadratic
programming (QP) framework for multi-agent optimization.
Our interest in quadratic programming is motivated
by problems in robotics~\cite{mellinger11} and 
data science~\cite{rodriguez10}, 
where some standard problems can
be formalized as QPs. 
The ``totally asynchronous'' label originates 
in~\cite{bertsekas1989parallel}, and it
describes a class of algorithms which tolerate arbitrarily
long delays, which our framework will do. In addition,
our developments will use block-based update laws in
which each agent updates only a small subset
of the decision variables in a problem, which reduces
each agent's computational burden and, as we will show,
reduces its onboard storage requirements as well. 
Both of these properties help enable the solution
of large-scale problems. 

Other work
on distributed quadratic programming 
includes~\cite{carli15,teixeira13,lee15,lee2015convergence,kozma2013distributed,todescato2015robust}.
Our work differs from these existing results because we
consider non-separable objective functions, and
because we consider
unstructured update laws (i.e., we do not require communications
and computations to occur in a particular sequence or pattern). 
Furthermore, we consider only deterministic problems, and our framework
converges exactly to a problem's solution, while some existing
works consider stochastic problems and converge approximately. 

Asynchrony in agents' communications and computations
implies that they will send and receive different information
at different times. As a result, they will disagree about
the values of decision variables in a problem. 
Just as it is difficult for agents to agree on this information,
it can also be difficult to agree on a stepsize value in
their algorithms. One could envision a network of agents solving
an agreement problem, e.g.,~\cite{ren05}, to compute a common 
stepsize,
though we instead allow agents to independently choose
stepsizes, subject to mild restrictions,
thereby eliminating the need to reach agreement
before optimizing. 

It has been shown that
regularizing problems can endow them with an inherent robustness
to asynchrony and improved convergence 
properties, e.g.,~\cite{koshal11,hale15,hale2017asynchronous}. 
Although regularizing is not required here,
we show, in a precise sense, that regularizing improves
convergence properties of our framework as well. 
It is common
for regularization-based approaches to require agents to
use the same regularization parameter, though this is undesirable
for the same reasons as using a common stepsize. Therefore,
we allow agents to independently choose regularization
parameters as well. 

To the best of our knowledge,
few works have considered both independent stepsizes
and regularizations. The most relevant is~\cite{koshal11},
which considers primal-dual algorithms for constrained
problems and synchronous primal updates. 
This paper is different in that we consider
unconstrained problems with totally asynchronous updates.
Regularizing of course introduces errors in a solution,
and we bound these errors in terms of agents' regularization
parameters. The end result is a totally asynchronous quadratic
programming framework in which communication delays can be
unbounded, and in which agents independently choose all stepsize
and regularization parameters. 

The rest of the paper is organized as follows. 
Section~\ref{sec:background} provides background on QPs
and formal problem statements.
Then, Section~\ref{sec:update} proposes an update law
to solve the problems of interest, and
Section~\ref{sec:convergence} proves its convergence.
Next, Section~\ref{sec:reg} shows that independent regularizations
lead to better-conditioned problems, and 
Section~\ref{sec:error} provides error bounds in terms of
agents' regularizations; Section~\ref{sec:discussion} provides
several corollaries and special cases for these bounds.  
After that, Section~\ref{sec:simulation}
provides simulation results, and Section~\ref{sec:conclusion}
concludes the paper.

\section{Background and Problem Statement} \label{sec:background}

In this section, we describe the quadratic optimization
problems to be solved, as well as the assumptions imposed upon these problems
and the agents that solve them. We then describe agents'
stepsizes and regularizations and introduce the need to allow agents
to choose these parameters independently. We next describe the benefits
of independent regularizations, and give two formal problem
statements that will be the focus of the remainder of the paper.

We consider a quadratic optimization problem distributed across a
network of $N$ agents, where agents are indexed over $i\in[N]:=\{1,...,N\}$.
Agent $i$ has a decision variable $x_{i}\in\mathbb{R}^{n_{i}},n_{i}\in\mathbb{N}$,
which we refer to as its state, and we allow for $n_{i}\neq n_{j}$
if $i\neq j$. The state $x_{i}$ is subject to the set constraint
$x_{i}\in X_{i}\subset\mathbb{R}^{n_{i}}$, and we make the following
assumption about each $X_{i}$.

\begin{assumption} \label{asm:setconst}
For all $i\in[N]$, the set $X_{i}\subset\mathbb{R}^{n_{i}}$
is non-empty, compact, and convex. $\hfill\triangle$
\end{assumption}

We define the network-level constraint set ${X:=X_{1}\times\cdots\times X_{N}}$,
and Assumption~\ref{asm:setconst} implies that $X$ is non-empty, compact, and convex.
We further define the global state as $x:=\left(x_{1}^{T},...,x_{N}^{T}\right)^{T}\in X\subset\mathbb{R}^{n}$,
where $n=\sum_{i\in[N]}n_{i}$. We consider quadratic objectives 
\begin{equation}
f(x):=\frac{1}{2}x^{T}Qx+r^{T}x,
\end{equation}
where $Q\in\mathbb{R}^{n\times n}$ and $r\in\mathbb{R}^{n}$. We
then make the following assumption about $f$.
\begin{assumption} \label{asm:Qpd}
In $f$, $Q=Q^{T}\succ0$. $\hfill\triangle$
\end{assumption}

Because $Q$ is positive definite, $f$ is strongly convex, and because
$f$ is quadratic, it is twice continuously differentiable, which
we indicate by writing that $f$ is $C^{2}$. In addition, ${\nabla f=Qx+r}$,
and~$\nabla f$ 
is therefore Lipschitz with constant $\|Q\|_{2}$. In this paper,
we divide $n\times n$ matrices into blocks. Given a matrix $B\in\mathbb{R}^{n\times n}$,
where $n=\sum_{i=1}^{N}n_{i}$, the $i^{th}$ block of $B$, denoted
$B^{[i]}$, is the $n_{i}\times n$ matrix formed by rows of $B$
with indices $\sum_{k=1}^{i-1}n_{k}+1$ through $\sum_{k=1}^{i}n_{k}$.
In other words, $B^{[1]}$ is the first $n_{1}$ rows of $B$, $B^{[2]}$
is the next $n_{2}$ rows, etc. Similarly, for a vector $b$, $b^{[1]}$
is the first $n_{1}$ entries of $b$, $b^{[2]}$ is the next $n_{2}$
entries, etc. Using this notion of a matrix block, we define $\nabla_{i}f:=\frac{\partial f}{\partial x_{i}}$,
and we see that $\nabla_{i}f(x)=Q^{[i]}x+r^{[i]}$.

Following our goal of reducing parametric coupling between agents,
we wish to allow agents to select stepsizes independently.
Allowing independent stepsizes will preclude the need for agents to
agree on a single value before optimizing, which gives agents a degree
of freedom in their update laws and eliminates the need to solve an
agreement problem before optimizing. Bearing this in mind, we state
the following problem, which will be one focus of the remainder of
the paper.

\textit{Problem 1:} Design a totally asynchronous distributed optimization
algorithm to solve
\begin{equation}
\underset{x\in X}{\text{minimize}}\quad\frac{1}{2}x^{T}Qx+r^{T}x,
\end{equation}
where only agent $i$ updates $x_{i}$, and where agents choose stepsizes
independently. $\hfill\triangle$

While an algorithm that satisfies the conditions stated in Problem
1 is sufficient to find a solution, there are additional characteristics
of $f$ that can be considered, namely $k_{Q}$, the spectral condition
number of $Q$. For a matrix $Q$ satisfying Assumption~\ref{asm:Qpd}, with eigenvalues
${\lambda_{1}(Q)\geq\lambda_{2}(Q)\geq...\geq\lambda_{n}(Q)}$, $k_{Q}$
is defined as ${k_{Q}:=\frac{\lambda_{1}(Q)}{\lambda_{n}(Q)}}$. In
standard centralized quadratic programming analyses, we find that $k_{Q}$
plays a vital role in determining convergence rates. In particular,
a problem with a large $k_{Q}$ is described as ``ill-conditioned,'' and
it will converge slower than a similar problem with a small $k_{Q}$.
We will find that this is indeed the case in Problem 1 as well.
Additionally, as will be shown below, $k_{Q}$ restricts agents' 
choices
of stepsizes. Therefore, a reduction in $k_{Q}$ will lead to a wider
range of stepsize choices for agents, in addition to faster convergence. 

Regularizations are commonly used for centralized quadratic programs
to improve $k_{Q}$ (i.e., to reduce it), and we will therefore use
them here. However, in keeping with the independence of agents' parameters,
we wish to allow agents to choose independent regularization parameters
as well. In particular, we should allow agent $i$ to use use the
regularization parameter $\alpha_{i}>0$, while allowing $\alpha_{i}\neq\alpha_{j}$
for $i\neq j$. The regularized form of $f$, denoted $f_{A}$,
is
\begin{equation}
f_{A}(x):=f(x)+\frac{1}{2}x^{T}Ax=\frac{1}{2}x^{T}(Q+A)x+r^{T}x,
\end{equation}
where $A=\text{diag}\left(\alpha_{1}I_{n_{1}},...,\alpha_{N}I_{n_{N}}\right)$,
and where $I_{n_{i}}$ is the $n_{i}\times n_{i}$ identity matrix.
Note that $\nabla_{i}f_{A}=Q^{[i]}x+r^{[i]}+\alpha_{i}x_{i}$, where 
we see that only $\alpha_{i}$
affects agent $i$'s updates.

With the goal of independent regularizations in mind, we now state
the second problem that we will solve.

\textit{Problem 2:} Design a totally asynchronous distributed optimization
algorithm to solve
\begin{equation}
\underset{x\in X}{\text{minimize}}\quad\frac{1}{2}x^{T}(Q+A)x+r^{T}x,
\end{equation}
where only agent $i$ updates $x_{i}$, and where agents independently choose their stepsizes
and regularizations. $\hfill\triangle$

Section III specifies the structure of the asynchronous communications
and computations used to solve Problem 1, and we will solve Problem
1 in Section IV. Afterwards, we will solve Problem 2 in Section V.

\section{Block-Based Multi-Agent Update Law} \label{sec:update}

To define the exact update law for each agent's state, we must first
describe the information stored onboard each agent and how agents
communicate with each other. Each agent will store a vector containing
its own state and that of every agent it communicates with. Formally,
we will denote agent $i$'s full vector of states by $x^{i}$, and
this is agent $i$'s local copy of the global state. Agent $i$'s
own states in this vector are denoted by $x_{i}^{i}$. The current
values stored onboard agent $i$ for agent $j$'s states are denoted
by $x_{j}^{i}.$ In the forthcoming update law, agent $i$ will only
compute updates for $x_{i}^{i}$, and it will share only $x_{i}^{i}$
with other agents when communicating. Agent $i$ will only change
the value of $x_{j}^{i}$ when agent $j$ sends its own state to agent
$i$. 

At time $k$, agent $i$\textquoteright s full state vector is denoted
$x^{i}(k)$, with its own states denoted $x_{i}^{i}(k)$ and those
of agent $j$ denoted $x_{j}^{i}(k)$. At any timestep, agent $i$
may or may not update its states due to asynchrony in agents\textquoteright{}
computations. As a result, we will in general have $x^{i}(k)\neq x^{j}(k)$
at all times $k$. We define the set $K^{i}$ to contain all times $k$ at which agent $i$ updates $x_{i}^{i}$; agent $i$
does not compute an update for time indices $k\notin K_{i}$. In designing
an update law, we must provide robustness to asynchrony while allowing
computations to be performed in a distributed fashion. First-order gradient 
descent methods are robust to many disturbances, and we therefore 
propose the following update law: 
\begin{equation}
x_{i}^{i}(k+1)=\begin{cases}
x_{i}^{i}(k)-\gamma_{i}\left(Q^{[i]}x^{i}(k)+r^{[i]}\right) & k\in K^{i}\\
x_{i}^{i}(k) & k\notin K^{i}
\end{cases},
\end{equation}
where agent $i$ uses some stepsize $\gamma_{i}>0$. This is equivalent
to agent $i$ using using the gradient descent law ${x_{i}^{i}(k+1)=x_{i}^{i}(k)-\gamma_{i}\nabla_{i}f\left(x^{i}(k)\right)}$
when it updates. The advantage of the block-based 
update law can be seen above, as agent $i$ only needs to know $Q^{[i]}$
and $r^{[i]}$. Requiring each agent to store the entirety of $Q$
and $r$ would require $O(n^{2})$ storage space, while $Q^{[i]}$
and $r^{[i]}$ only require $O(n)$. For large quadratic programs,
this block-based update law dramatically reduces each agent's onboard
storage requirements.

In order to account for communication delays, we use $\tau_{j}^{i}(k)$
to denote the time at which the value of $x_{j}^{i}(k)$ was originally
computed by agent $j$. For example, if agent $j$ computes a state
update at time $k_{a}$ and immediately transmits it to agent $i$,
then agent $i$ may receive this state update at time $k_{b}>k_{a}$
due to communication delays. Then $\tau_{j}^{i}$ is defined so that
$\tau_{j}^{i}(k_{b})=k_{a}$, which relates the time of receipt by
agent $i$ to the time at which agent $j$ originally computed the
piece of data being sent. For $K^{i}$ and $\tau_{j}^{i}$, we assume
the following.

\begin{assumption} \label{asm:infupdate}
For all $i\in[N]$, the set $K^{i}$ is infinite.
Moreover, for all $i\in[N]$ and $j\in[N]\backslash\{i\}$, if $\left\{ k_{d}\right\} _{d\in\mathbb{N}}$
is a sequence in $K^{i}$ tending to infinity, then
\end{assumption}
\begin{equation}
\lim_{d\rightarrow\infty}\tau_{j}^{i}(k_{d})=\infty. \tag*{$\triangle$}
\end{equation}
Assumption~\ref{asm:infupdate} is simply a formalization of the requirement that no
agent ever permanently stop updating and sharing its own state with
any other agent. For $i\neq j$, the sets $K^{i}$ and $K^{j}$ do
not need to have any relationship because agents' updates are asynchronous.
Our proposed update law for all agents can then be written as follows.

\textit{Algorithm 1:} For all $i\in[N]$ and $j\in[N]\backslash\{i\}$,
execute
\begin{align*}
x_{i}^{i}(k+1) & =\begin{cases}
x_{i}^{i}(k)-\gamma_{i}\left(Q^{[i]}x^{i}(k)+r^{[i]}\right) & k\in K^{i}\\
x_{i}^{i}(k) & k\notin K^{i}
\end{cases}\\
x_{j}^{i}(k+1) & =\begin{cases}
x_{j}^{j}\left(\tau_{j}^{i}(k+1)\right) & \text{i receives j's state at k+1}\\
x_{j}^{i}(k) & \text{otherwise} \hfill\diamond
\end{cases}
\end{align*}

In Algorithm 1 we see that $x_{j}^{i}$ changes only when agent $i$
receives a transmission directly from agent $j$; otherwise it remains
constant. This implies that agent $i$ can update its own state using
an old value of agent $j$\textquoteright s state multiple times and
can reuse different agents\textquoteright{} states different numbers
of times. Showing convergence of this update law must 
account for these delays, in addition to providing  stepsize
bounds for each agent, and that is the subject of the next section.

\section{Convergence of Asynchronous Optimization}\label{sec:convergence}

In this section, we prove convergence of the multi-agent block update
law in Algorithm 1. This will be shown using a block-maximum norm 
to measure convergence, along with 
a collection of nested sets to show Lyapunov-like convergence.
We will derive stepsize bounds from these concepts that will be used
to show asymptotic convergence of all agents.

\subsection{Block-Maximum Norms}

The convergence of Algorithm 1 will be measured using a block-maximum
norm as 
in~\cite{bertsekas1989convergence},~\cite{bertsekas1989parallel}, 
and~\cite{hale2017asynchronous}. We do this to permit agents
to use independent normalizations in Problem 1 to weight different
components of $x$ differently when estimating convergence to an optimum. Below, we refer to~$x_{i}^{i}$ as the~$i^{th}$ block of~$x^{i}$ and~$x_{j}^{i}$
as the~$j^{th}$ block of~$x^{i}$. We next define the block-maximum
norm that will be used to measure convergence.

\textit{Definition 1:} Let $x\in\mathbb{R}^{n}$ consist of $N$ blocks,
with ${x_{i}\in\mathbb{R}^{n_{i}}}$ being the $i^{th}$ block. The
$i^{th}$ block is weighted by some normalization constant $\omega_{i}\geq1$
and is measured in the $p_{i}$-norm for some $p_{i}\in[1,\infty]$.
The norm of the full vector $x$ is defined as the maximum norm of
any single block, i.e.,
\begin{equation}
\|x\|_{max}:=\max_{i\in[N]}\frac{\|x_{i}\|_{p_{i}}}{\omega_{i}}. \tag*{$\blacktriangle$}
\end{equation}

The following lemma allows us to upper-bound the induced block-maximum
matrix norm by the Euclidean matrix norm, which will be used below
in our convergence analysis. 

\begin{lemma} \label{lem:blockmax}
Suppose for all $i\in[N]$ that agent $i$ uses the
weight $\omega_{i} \geq 1$ and $p_{i}$-norm,~$p_i \in [1, \infty]$, 
in the above block-maximum
norm. Let $p_{min}:=\min_{i\in[N]}p_{i}$ and let $\omega_{min}:=\min_{i\in[N]}\omega_{i}$.
Then for all $B\in\mathbb{R}^{n\times n}$,
\end{lemma}
\begin{equation}
\|B\|_{max}\leq\begin{cases}
n^{\left(p_{min}^{-1}-\frac{1}{2}\right)}\omega_{min}^{-1}\|B\|_{2} & p_{min}<2\\
\omega_{min}^{-1}\|B\|_{2} & p_{min}\geq2
\end{cases}.
\end{equation}

\textit{Proof:} See Lemma 1 in \cite{hochhaus2018asynchronous}. \hfill$\blacksquare$
\subsection{Convergence Via Lyapunov Sub-Level Sets}

We will now analyze the convergence of Algorithm 1 when agents are
communicating asynchronously. In order to show convergence, we construct
a sequence of sets, $\{X(s)\}_{s\in\mathbb{N}}$, based on work 
in~\cite{bertsekas1989convergence} 
and~\cite{bertsekas1989parallel}. These sets behave analogously to sub-level sets
of a Lyapunov function, and they will enable an invariance type argument
in our convergence proof. Below, we use $\hat{x}:=\arg\min_{x\in X}f(x)$
for the minimizer of $f$. 
For simplicity, we assume that~$\hat{x}$ is in the interior of~$X$,
though all of our results hold without modification if a 
projection onto~$X_i$ is added to agent~$i$'s update
in Algorithm~$1$, and that is the only change required
if~$\hat{x}$ is not in the interior of~$X$. 
We state the following assumption on these
sets, and below we will construct a sequence of sets that satisfies
this assumption.

\begin{assumption} \label{asm:setsexist}
There exists a collection of sets $\{X(s)\}_{s\in\mathbb{N}}$
that satisfies:

1) $...\subset X(s+1)\subset X(s)\subset...\subset X$

2) $\lim_{s\rightarrow\infty}X(s)=\left\{ \hat{x}\right\} $

3) There exists $X_{i}(s)\subset X_{i}$ for all $i\in[N]$ and $s\in\mathbb{N}$
such that $X(s)=X_{1}(s)\times...\times X_{N}(s)$

4) $\theta_{i}(y)\in X_{i}(s+1)$, where $\theta_{i}(y):=y_{i}-\gamma_{i}\nabla_{i}f(y)$
for all $y\in X(s)$ and $i\in[N]$. $\hfill\triangle$
\end{assumption}

Assumptions~\ref{asm:setsexist}.1 and~\ref{asm:setsexist}.2 jointly guarantee that the collection $\{X(s)\}_{s\in\mathbb{N}}$
is nested and that they converge to a singleton containing $\hat{x}$.
Assumption~\ref{asm:setsexist}.3 allows for the blocks 
of~$x$ to be updated independently by
the agents, which allows for decoupled update laws. 
Assumption~\ref{asm:setsexist}.4
ensures that state updates make only forward progress toward $\hat{x}$,
which ensures that each set is forward-invariant in time. It is shown
in~\cite{bertsekas1989convergence} and~\cite{bertsekas1989parallel} that the existence of such a sequence of sets
implies asymptotic convergence of the asynchronous update law in Algorithm
1. We therefore use this strategy to show asymptotic convergence
in this paper. We propose to use the construction 
\begin{equation} \label{eqn:setcon}
X(s)=\left\{ y\in X:\|y-\hat{x}\|_{max}\leq q^{s}nD_{o}\right\} ,
\end{equation}
where we define~$D_{o}:=\max_{i\in[N]}\|x^{i}(0)-\hat{x}\|_{max}$,
which is the block furthest from $\hat{x}$ onboard any agent at
timestep zero, and where we define the constant
\begin{equation}
q=\|I-\Gamma Q\|_{2},
\end{equation}
with $\Gamma=\text{diag}\left(\gamma_{1}I_{n_{1}},...,\gamma_{N}I_{n_{N}}\right)$.
We will use the fact that each update contracts towards $\hat{x}$ by a factor
of $q$, and the following theorem will establish bounds on every $\gamma_{i}$ 
that imply $q\in(0,1)$. This result will be used to show
convergence of Algorithm 1 through satisfaction of 
Assumption~\ref{asm:setsexist}.

\begin{theorem} \label{thm:step}
Let $Q=Q^{T}\succ0$,~$Q \in \R^{n \times n}$, have condition number $k_{Q}$, and let $\Gamma=\text{diag}\left(\gamma_{1}I_{n_{1}},...,\gamma_{N}I_{n_{N}}\right)$.
If
\begin{equation} \label{eqn:stepbound}
\gamma_{i}\in\left(\frac{\sqrt{k_{Q}}-1}{||Q||_{2}\sqrt{k_{Q}}},\frac{\sqrt{k_{Q}}+1}{||Q||_{2}\sqrt{k_{Q}}}\right)\,\,\,\text{ for all  }\,\,\,i\in[N],
\end{equation}
then $||I-\Gamma Q||_{2}<1$.
\end{theorem}

\textit{Proof:} To avoid interrupting the flow of the paper, proof
of Theorem~\ref{thm:step} can be found in the appendix.$\hfill\blacksquare$

In Theorem~\ref{thm:step}, we note that any choice of $\gamma_{lower}$ and $\gamma_{upper}$
that satisfies 
\begin{equation} \label{eqn:yupper}
\begin{split} 
\hspace{-0.3cm}\frac{1}{2\gamma_{upper}\|Q\|_{2}}\hspace{-0.05cm}\left(\hspace{-0.05cm}\frac{(\sqrt{k_{Q}}\hspace{-0.05cm}+\hspace{-0.05cm}1)^{2}}{\sqrt{k_{Q}}}
\hspace{-0.05cm}-\hspace{-0.05cm}\frac{\gamma_{upper}}{\gamma_{lower}}\frac{(\sqrt{k_{Q}}\hspace{-0.05cm}-\hspace{-0.05cm}1)^{2}}{\sqrt{k_{Q}}}\hspace{-0.05cm}\right)\hspace{-0.05cm}>\hspace{-0.05cm}1
\end{split}
\end{equation}
defines a valid interval of step sizes to ensure $\|I-\Gamma Q\|_{2}<1$.
While the algebra is omitted here, it can be shown that the range
defined in Theorem~\ref{thm:step} is the largest of these intervals. Note also
that due to the structure of Equation~\eqref{eqn:yupper}, if the exact values of
$k_{Q}$ and $\|Q\|_{2}$ are unavailable or difficult to compute,
then they can be replaced in Equation~\eqref{eqn:stepbound} by upper bounds. These could
include using Gershgorin's Circle Theorem or the trace of $Q$ to
bound $\|Q\|_{2}$ and using results such as in ~\cite{cheng2014note} to bound $k_{Q}$.

Letting $\gamma_{i}\in\left(\frac{\sqrt{k_{Q}}-1}{\|Q\|_{2}\sqrt{k_{Q}}},\frac{\sqrt{k_{Q}}+1}{\|Q\|_{2}\sqrt{k_{Q}}}\right)$
for all $i\in[N]$ and recalling our construction of sets $\left\{ X(s)\right\} _{s\in\mathbb{N}}$
as
\begin{equation}
X(s)=\left\{ y\in X:\|y-\hat{x}\|_{max}\leq q^{s}nD_{o}\right\} ,
\end{equation}
we next show that Assumption~\ref{asm:setsexist} is satisfied, thereby ensuring
convergence of Algorithm 1.

\begin{theorem} \label{thm:setswork}
If $\Gamma$ satisfies the conditions in Theorem
1, then the collection of sets $\left\{ X(s)\right\} _{s\in\mathbb{N}}$
as defined in Equation~\eqref{eqn:setcon} satisfies Assumption~\ref{asm:setsexist}.
\end{theorem}
\textit{Proof:} For Assumption~\ref{asm:setsexist}.1, by definition we have
\begin{equation}
X(s+1)=\left\{ y\in X:\|y-\hat{x}\|_{max}\leq q^{s+1}nD_{o}\right\} .
\end{equation}
Since $q\in(0,1)$, we have $q^{s+1}<q^{s}$, which results in~${\|y-\hat{x}\|_{max}\leq q^{s+1}nD_{o}<q^{s}nD_{o}}$.
Then $y\in X(s+1)$ implies $y\in X(s)$ and $X(s+1)\subset X(s)\subset X$,
as desired.

For Assumption~\ref{asm:setsexist}.2 we find
\begin{equation}
\lim_{s\rightarrow\infty}X(s)=\lim_{s\rightarrow\infty}\left\{ y\in X:\|y-\hat{x}\|_{max}\leq q^{s}nD_{o}\right\} =\left\{ \hat{x}\right\}.
\end{equation}
The structure of the weighted block-maximum
norm then allows us to see that $\|y-\hat{x}\|_{max}\leq q^{s}nD_{o}$
if and only if $\frac{1}{\omega_{i}}\|y_{i}-\hat{x}_{i}\|_{p_{i}}\leq q^{s}nD_{o}$
for all $i\in[N].$ It then follows that
\begin{equation}
X_{i}(s)=\left\{ y_{i}\in X_{i}:\frac{1}{\omega_{i}}\|y_{i}-\hat{x}_{i}\|_{p_{i}}\leq q^{s}nD_{o}\right\} ,
\end{equation}
which gives $X(s)=X_{1}(s)\times...\times X_{N}(s)$, thus satisfying
Assumption~\ref{asm:setsexist}.3.

We then see that, for $y\in X(s)$,
\begin{equation}
\frac{\|\theta_{i}(y)-\hat{x}_{i}\|_{p_{i}}}{\omega_{i}}= \frac{1}{\omega_{i}}\Big\|y_{i}- \gamma_{i}\left(Q^{[i]}y+r^{[i]}\right)-\hat{x}_{i}+\gamma_{i}\left(Q^{[i]}\hat{x}+r^{[i]}\right)\Big\|_{p_{i}},
\end{equation}
which follows from the definition of $\theta_{i}(y)$ and the fact
that $\nabla_{i}f(\hat{x})=0$. We then find
\begin{align*}
\frac{\|\theta_{i}(y)-\hat{x}_{i}\|_{p_{i}}}{\omega_{i}} & \leq\max_{i\in[N]}\frac{1}{\omega_{i}}\|y_{i}-\gamma_{i} \left(Q^{[i]}y+r^{[i]}\right)-\hat{x}_{i}+\gamma_{i}\left(Q^{[i]}\hat{x}+r^{[i]}\right)\|_{p_{i}}\\
& =\|y-\hat{x}-\Gamma\left(Qy+r\right)+\Gamma\left(Q\hat{x}+r\right)\|_{max},
\end{align*}
which follows from our definition of the block-maximum norm. Continuing, we find
\begin{align*}
\frac{\|\theta_{i}(y)-\hat{x}_{i}\|_{p_{i}}}{\omega_{i}} 
 &\leq\|y-\hat{x}-\Gamma Q\left(y-\hat{x}\right)\|_{max}\\
 & \leq\|I-\Gamma Q\|_{max}\|y-\hat{x}\|_{max}\\
 &\leq\begin{cases}
\frac{n^{\left(p_{min}^{-1}-\frac{1}{2}\right)}}{\omega_{min}}\|I-\Gamma Q\|_{2}\|y-\hat{x}\|_{max} & p_{min}<2\\
\frac{1}{\omega_{min}}\|I-\Gamma Q\|_{2}\|y-\hat{x}\|_{max} & p_{min}\geq2
\end{cases},
\end{align*}
where the last inequality follows from Lemma~\ref{lem:blockmax}. Seeing that $\|I-\Gamma Q\|_{2}=q\in(0,1)$,
and using the hypothesis that ${y\in X(s)}$, we find
\begin{align*}
\frac{\|\theta_{i}(y)-\hat{x}_{i}\|_{p_{i}}}{\omega_{i}} & \leq\begin{cases}
\frac{n^{\left(p_{min}^{-1}-\frac{1}{2}\right)}}{\omega_{min}}q\|y-\hat{x}\|_{max} & p_{min}<2\\
\frac{1}{\omega_{min}}q\|y-\hat{x}\|_{max} & p_{min}\geq2
\end{cases}\\
 & \leq\begin{cases}
\frac{n^{\left(p_{min}^{-1}-\frac{1}{2}\right)}}{\omega_{min}}q^{s+1}D_{o} & p_{min}<2\\
\frac{1}{\omega_{min}}q^{s+1}D_{o} & p_{min}\geq2
\end{cases}\\
 & \leq\begin{cases}
q^{s+1}nD_{o} & p_{min}<2\\
q^{s+1}nD_{o} & p_{min}\geq2
\end{cases},
\end{align*}
where the bottom case follows from $\omega_{min}\geq1$ and the top
case follows from $\omega_{min}\geq1$ 
and~$p_{min}^{-1}-\frac{1}{2}<1$
(since $p_{i}\in[1,\infty]$ for all~$i \in [N]$). 
Then $\theta_{i}(y)\in X_{i}(s+1)$ and
Assumption~\ref{asm:setsexist}.4 is satisfied. $\hfill\blacksquare$

Regarding Problem 1, we therefore state the following:

\begin{theorem} \label{thm:alg1works}
Algorithm 1 solves Problem 1 and asymptotically
converges to $\hat{x}$.
\end{theorem}

\textit{Proof:} Theorem~\ref{thm:setswork} shows the construction of the sets $\left\{ X(s)\right\} _{s\in\mathbb{N}}$
satisfies Assumption~\ref{asm:setsexist}, and from ~\cite{bertsekas1989convergence} and ~\cite{bertsekas1989parallel} we see this implies
convergence of Algorithm 1 for all $i\in[N]$. The total
asynchrony required by Problem 1 is incorporated by not requiring
delay bounds, and agents do not require any coordination in selecting
stepsizes, which means that all of the criteria of Problem 1 are satisfied.
$\hfill\blacksquare$

From these requirements, we see that agent $i$ only needs to be
initialized with $Q^{[i]}$, $r^{[i]}$, and upper bounds on $\|Q\|_{2}$
and $k_{Q}$, which can be set by a network operator. Agents are then
free to choose normalizations and stepsizes independently, provided
stepsizes obey the bounds established in Theorem~\ref{thm:step}. 

\section{Independently Regularized Quadratic Programs} \label{sec:reg}
Equation~\eqref{eqn:stepbound} quantifies the relationship 
between the condition
number of~$Q$ and 
agents' stepsize bounds. For a perfectly conditioned
matrix $Q$, such as the identity, we have $k_{Q}=1$, which implies that
stepsizes may be chosen from the open interval $\left(0,\frac{2}{\|Q\|_{2}}\right)$.
This is the familiar $\frac{2}{L}$ bound encountered in conventional gradient descent
settings. As the problem becomes more ill-conditioned and $k_{Q}$
increases, the allowable interval of stepsize choices becomes narrower
in Equation~\eqref{eqn:stepbound}. Specifically, as $k_{Q}\rightarrow\infty$, the set
of allowable stepsizes approaches the degenerate interval $\left[\frac{1}{\|Q\|_{2}},\frac{1}{\|Q\|_{2}}\right]$,
implying that all stepsizes must be equal in such a case. 

Equation~\eqref{eqn:stepbound} suggests that if the condition number of a quadratic program is
reduced, then the interval of allowable step sizes can be lengthened.
In particular, regularizing $f$ can improve $k_{Q}$ to do so. 
As stated in
Problem 2, we want to allow agents to choose regularization parameters
independently. Before we analyze the effects of independently chosen
regularizations on $k_{Q}$, we must first show that an algorithm
that utilizes them will preserve the convergence properties of Algorithm~1. As shown above, a regularized cost function takes the form
\begin{equation}
f_{A}(x):=\frac{1}{2}x^{T}(Q+A)x+r^{T}x,
\end{equation}
where $Q+A$ is symmetric positive definite because ${Q=Q^{T}\succ0}$.
We now state the following theorem that confirms that minimizing $f_{A}$
succeeds.

\begin{theorem} \label{thm:prob2solved}
Suppose that $A\succ0$ is diagonal with positive diagonal
entries. 
Then Algorithm
1 satisfies the conditions stated in Problem 2 when $f_{A}$ is minimized.
\end{theorem}

\textit{Proof:} Replacing $Q$ with $Q+A,$ all assumptions and conditions
used to prove Theorem~\ref{thm:alg1works} hold, with the only modifications being that
the network will converge to ${\hat{x}_{A}:=\arg\min_{x\in X}f_{A}(x)}$,
and the agents must now be initialized with $Q^{[i]}$, $r^{[i]}$,
and upper limits on $\|Q+A\|_{2}$ and the condition number of $Q+A$.
$\hfill\blacksquare$

Theorem~\ref{thm:prob2solved} implies that there must be some known upper bound on $\|Q+A\|_{2}$
and the condition number of $Q+A$. These bounds can be determined using
bounds on agents' allowable regularization parameters, which we develop 
now. First we will establish the
following theorem, which demonstrates that independent regularizations
can indeed reduce the condition number of the quadratic program.

\begin{theorem} \label{thm:kbound}
Let there be two $n\times n$ matrices ${Q=Q^{T}\succ0}$
and $A=\text{diag}\left(\alpha_{1}I_{n_{1}},...,\alpha_{N}I_{n_{N}}\right)\succ0$
with respective condition numbers $k_{Q}=\frac{\lambda_{1}(Q)}{\lambda_{n}(Q)}$
and $k_{A}=\frac{\lambda_{1}(A)}{\lambda_{n}(A)}=\frac{\alpha_{max}}{\alpha_{min}}$,
where $\alpha_{max}=\max_{i}\alpha_{i}$ and $\alpha_{min}=\min_{i}\alpha_{i}$.
If $\frac{\alpha_{max}}{\alpha_{min}}<k_{Q}$, then $k_{Q+A}<k_{Q}.$
\end{theorem}

\textit{Proof:} Using $\frac{\alpha_{max}}{\alpha_{min}} = k_A <k_{Q}$,
we find 
$\frac{\lambda_{1}(A)}{\lambda_{n}(A)} <\frac{\lambda_{1}(Q)}{\lambda_{n}(Q)}$.
Rearranging, we find~${\lambda_{n}(Q)\lambda_{1}(A)<\lambda_{1}(Q)\lambda_{n}(A)}$.
Adding $\lambda_{1}(Q)\lambda_{n}(Q)$ to both sides and factoring gives
\begin{equation}
\lambda_{n}(Q)\left(\lambda_{1}(Q)+\lambda_{1}(A)\right)<\lambda_{1}(Q)\left(\lambda_{n}(Q)+\lambda_{n}(A)\right),
\end{equation}
which we rearrange again to find
\begin{equation} \label{eqn:weyllambda}
\frac{\lambda_{1}(Q)+\lambda_{1}(A)}{\lambda_{n}(Q)+\lambda_{n}(A)}<\frac{\lambda_{1}(Q)}{\lambda_{n}(Q)}.
\end{equation}
From Weyl's inequalities, if $B$ and $C$ are $n\times n$ Hermitian
matrices, then $\lambda_{1}(B+C)\leq\lambda_{1}(B)+\lambda_{1}(C)$ and
${\lambda_{n}(B)+\lambda_{n}(C)\leq\lambda_{n}(B+C)}$ ~\cite[Fact 5.12.2]{bernstein2009matrix}. Therefore
\begin{equation} 
\frac{\lambda_{1}(Q+A)}{\lambda_{n}(Q+A)} \leq\frac{\lambda_{1}(Q)+\lambda_{1}(A)}{\lambda_{n}(Q)+\lambda_{n}(A)}.
\end{equation}
Combining this with Equation~\eqref{eqn:weyllambda}
completes the proof. 
\hfill $\blacksquare$

Since the regularization matrix $A$ is chosen by the agents, it is
always possible to have $k_{A}<k_{Q}$ provided $k_{Q}>1$. Of course,
we will only regularize such a problem, and thus, for all practical
purposes, regularizing improves the conditioning of our problems as
long as $k_{A}$ is better-conditioned than $k_{Q}$. 
If we wish to reduce $k_{Q+A}$ such that it is bounded above by some
desired condition number $k_{D}$, we can bound 
$\alpha_{min}$ via 
\begin{equation} \label{eqn:aminbound}
\hspace{-0.3cm}\alpha_{min}\hspace{-0.05cm}>\hspace{-0.05cm}\|Q\|_{2}\hspace{-0.05cm}\left(k_{D}^{-1}\hspace{-0.05cm}-\hspace{-0.05cm}k_{Q}^{-1}\hspace{-0.05cm}\right)\hspace{-0.05cm}+\frac{\epsilon\|Q\|_{2}^{2}}{k_{Q}k_{D}\left(\|r\|_{2}k_{Q}\hspace{-0.05cm}-\hspace{-0.05cm}\epsilon\|Q\|_{2}\right)},
\end{equation}
which works for any $k_{D}$ such that
\begin{equation} \label{eqn:kdbound}
k_{D}>k_{Q}-\frac{\epsilon\|Q\|_{2}\left(k_{Q}-1\right)}{\|r\|_{2}k_{Q}}.
\end{equation}

Theorem~\ref{thm:kbound} suggests
that if we want to expand the interval of allowable stepsizes via
regularizing, then we should choose large, homogeneous regularization
parameters among agents. However, larger regularizations will lead
to larger errors in the solution to a quadratic program, and the ability
to take larger steps must be balanced with the quality of solution
obtained. In the next section, we will quantify the relationship between
regularizations and error, and develop a set of guidelines to govern
agents' selection of regularization parameters based
on desired error bounds. 

\section{Regularization Error Bound} \label{sec:error}

We see from the structure of the quadratic program that the solution
to the unregularized minimization problem is~${\hat{x}=-Q^{-1}r}$,
where $Q^{-1}$ is well-defined because ${Q\succ0}$. Similarly, the solution to the regularized minimization problem 
is~$\hat{x}_{A}=-(Q+A)^{-1}r$.
We therefore define the regularization error when using regularization matrix
$A$ as~$e_{A}:=\|\hat{x}-\hat{x}_{A}\|_{2}$.
We now upper bound $e_{A}$ in terms of the entries of~$A$.

\begin{theorem} \label{thm:regerr}
For the agent-specified regularization matrix
$A=\text{diag}\left(\alpha_{1}I_{n_{1}},...,\alpha_{N}I_{n_{N}}\right)$,
the regularization error ${e_{A}:=\|\hat{x}-\hat{x}_{A}\|_{2}}$ is
bounded via 
\begin{equation} \label{eqn:errbound}
e_{A}\leq\frac{\|r\|_{2}k_{Q}^{2}\alpha_{max}}{\|Q\|_{2}^{2}+\|Q\|_{2}k_{Q}\alpha_{max}},
\end{equation}
where $\alpha_{max}:=\max_{i}\alpha_{i}$.
\end{theorem}

\textit{Proof:} By definition, we have
\begin{equation}
e_{A} =\|Q^{-1}r-(Q+A)^{-1}r\|_{2} \leq\|Q^{-1}-(Q+A)^{-1}\|_{2}\|r\|_{2}.
\end{equation}
The above inequality can be rewritten using Lemma~\ref{lem:woodbury}.

\begin{lemma} \label{lem:woodbury}
Let $Q$ and $A$ be positive definite matrices. Then
\begin{equation}
Q^{-1}-(Q+A)^{-1}=(QA^{-1}Q+Q)^{-1}.
\end{equation}
\end{lemma}

\textit{Proof:} From the Woodbury matrix identity,
\begin{equation}
(B+UCV)^{-1}=B^{-1}-B^{-1}U(C^{-1}+VB^{-1}U)^{-1}VB^{-1}.
\end{equation}
Let $B=U=C=Q$ and $C=Q^{-1}AQ^{-1}.$ Then
\begin{equation}
(Q+A)^{-1}=Q^{-1}-Q^{-1}Q(QA^{-1}Q+QQ^{-1}Q)^{-1}QQ^{-1},
\end{equation}
which simplifies to~$(Q+A)^{-1}=Q^{-1}-(QA^{-1}Q+Q)^{-1}$.
Rearranging completes the lemma. \hfill $\blacksquare$

Therefore,~${e_{A}\leq\|(QA^{-1}Q+Q)^{-1}\|_{2}\|r\|_{2}}$ using Lemma~\ref{lem:woodbury}.
Note that $QA^{-1}Q+Q$ is a symmetric positive definite matrix. For
any symmetric positive definite matrix $M$, 
we know~$\|M^{-1}\|_{2}=\lambda_{1}(M^{-1})=\lambda_{n}^{-1}(M)$.
Applying this to~$(QA^{-1}Q+Q)^{-1}$
gives
\begin{equation}
e_{A}\leq\frac{\|r\|_{2}}{\lambda_{n}(QA^{-1}Q+Q)}.
\end{equation}
From Weyl's Inequalities \cite[Fact 5.12.2]{bernstein2009matrix} we then have
\begin{equation}
e_{A}\leq\frac{\|r\|_{2}}{\lambda_{n}(QA^{-1}Q)+\lambda_{n}(Q)}.
\end{equation}
Since $QA^{-1}Q$ is similar to $Q^{2}A^{-1}$, we know that
$\lambda_{n}(QA^{-1}Q)=\lambda_{n}(Q^{2}A^{-1})$,
and since $Q^{2}$ and $A^{-1}$ are both symmetric positive definite,
and $\lambda_{n}(Q^{2})=\lambda_{n}^{2}(Q)$, we can say ~\cite[Fact 8.19.18]{bernstein2009matrix} 
that~$\lambda_{n}^{2}(Q)\lambda_{n}(A^{-1})\leq\lambda_{n}(Q^{2}A^{-1})$,
which gives
\begin{equation}
e_{A} \leq\frac{\|r\|_{2}}{\lambda_{n}^{2}(Q)\lambda_{n}(A^{-1})+\lambda_{n}(Q)}  =\frac{\|r\|_{2}}{\lambda_{n}^{2}(Q)\alpha_{max}^{-1}+\lambda_{n}(Q)}=\frac{\|r\|_{2}\alpha_{max}}{\lambda_{n}^{2}(Q)+\lambda_{n}(Q)\alpha_{max}}.
\end{equation}
Substituting in $\lambda_{n}(Q)=\frac{\|Q\|_{2}}{k_{Q}}$ completes the proof. \hfill $\blacksquare$

\section{Discussion of Results} \label{sec:discussion}

In this section, we briefly provide some remarks and interpretation
of our results. We begin with the following corollary that follows
from Theorem~\ref{thm:regerr}.
\begin{corollary}
As $\alpha_{max}\rightarrow0$, $e_{A}\rightarrow0$.
\end{corollary}
It is desirable for an error bound to be sharp in the sense that equality
is achieved for at least one point. Corollary~1 shows that this is
achieved in Equation~\eqref{eqn:errbound}. This implies zero regularization error when
zero regularization is applied, which establishes 
Equation~\eqref{eqn:errbound} as a
sharp upper bound.
\begin{corollary}
As $\alpha_{max}\rightarrow\infty$, $e_{A}\rightarrow\frac{\|r\|_{2}k_{Q}}{\|Q\|_{2}}$.
\end{corollary}
Another notable behavior of Equation~\eqref{eqn:errbound} is a finite upper bound on
the regularization error, as shown in Corollary~2. To elaborate on
this, consider the case where ${A=\alpha I}$, where $\alpha$ is a
positive scalar. As $\alpha\rightarrow\infty$, we see ${\hat{x}_{A}=-(Q+A)^{-1}r\rightarrow0}$.
That is, very large regularizations move the regularized solution
$\hat{x}_{A}$ towards zero. We also see as $\alpha\rightarrow\infty,$
$e_{A}\rightarrow\|Q^{-1}r\|_{2},$ since the norm of the error between
zero and the unregularized solution is simply the norm of the unregularized
solution itself. The upper bound from Equation~\eqref{eqn:errbound} as $\alpha_{max}\rightarrow\infty$
is due to the fact that
\begin{equation}
\|Q^{-1}r\|_{2}\leq\|Q^{-1}\|_{2}\|r\|_{2}=\frac{\|r\|_{2}k_{Q}}{\|Q\|_{2}}.
\end{equation}
If there is some desired upper bound $\epsilon$ on the error $e_{A}$,
then an upper bound on $\alpha_{max}$ can be determined as
\begin{equation} \label{eqn:amaxbound}
\alpha_{max}<\frac{\epsilon\|Q\|_{2}^{2}}{\|r\|_{2}k_{Q}^{2}-\epsilon\|Q\|_{2}k_{Q}},
\end{equation}
which holds as long as $\epsilon<\frac{\|r\|_{2}k_{Q}}{\|Q\|_{2}}$.
This bound can be combined with the bound on~$\alpha_{min}$
in Equation~\eqref{eqn:aminbound} to present agents
with an admissible range of regularization parameters, which
we do in the next section. 

\section{Simulation} \label{sec:simulation}
In this section, two simulations are run in which a 
network of agents solves
a quadratic program using Algorithm 1. In the first simulation, the
problem is unregularized and agents select stepsizes independently
from the allowable set for the problem, and they 
then asynchronously solve
the quadratic program. In the second simulation, 
a desired condition number bound
and regularization error bound are specified.
Agents independently select regularization
parameters to satisfy those requirements. Agents then choose stepsizes
independently from the expanded allowable set, and then asynchronously
solve the quadratic program. The desired condition number bound and
error bound conditions are shown to be satisfied, and the convergence
histories of the first and second simulations are compared.

The network consists of 25 agents, each with 4 states.
The quadratic program is randomly generated with $k_{Q}=100$ and $\|Q\|_{2}=100$.
For the first simulation, agents randomly choose stepsizes from the
interval $(0.009,0.0110)$, which is obtained from Equation ~\eqref{eqn:stepbound} and
the values of $k_{Q}$ and $\|Q\|_{2}$. In the second, we specify
the desired condition number ${k_{D}=10}$ and the upper bound
$\epsilon=0.1$ on regularization error. 
Agents then randomly choose regularization parameters
from the interval $(11,20)$, which follows from 
Equations~\eqref{eqn:amaxbound} and~\eqref{eqn:aminbound}. Agents then randomly select stepsizes from the
interval $(0.0056,0.0117)$, which is obtained as before.

\setlength{\belowcaptionskip}{-40pt}
\begin{figure}[!tp]
\centering
\includegraphics[width=3.5in]{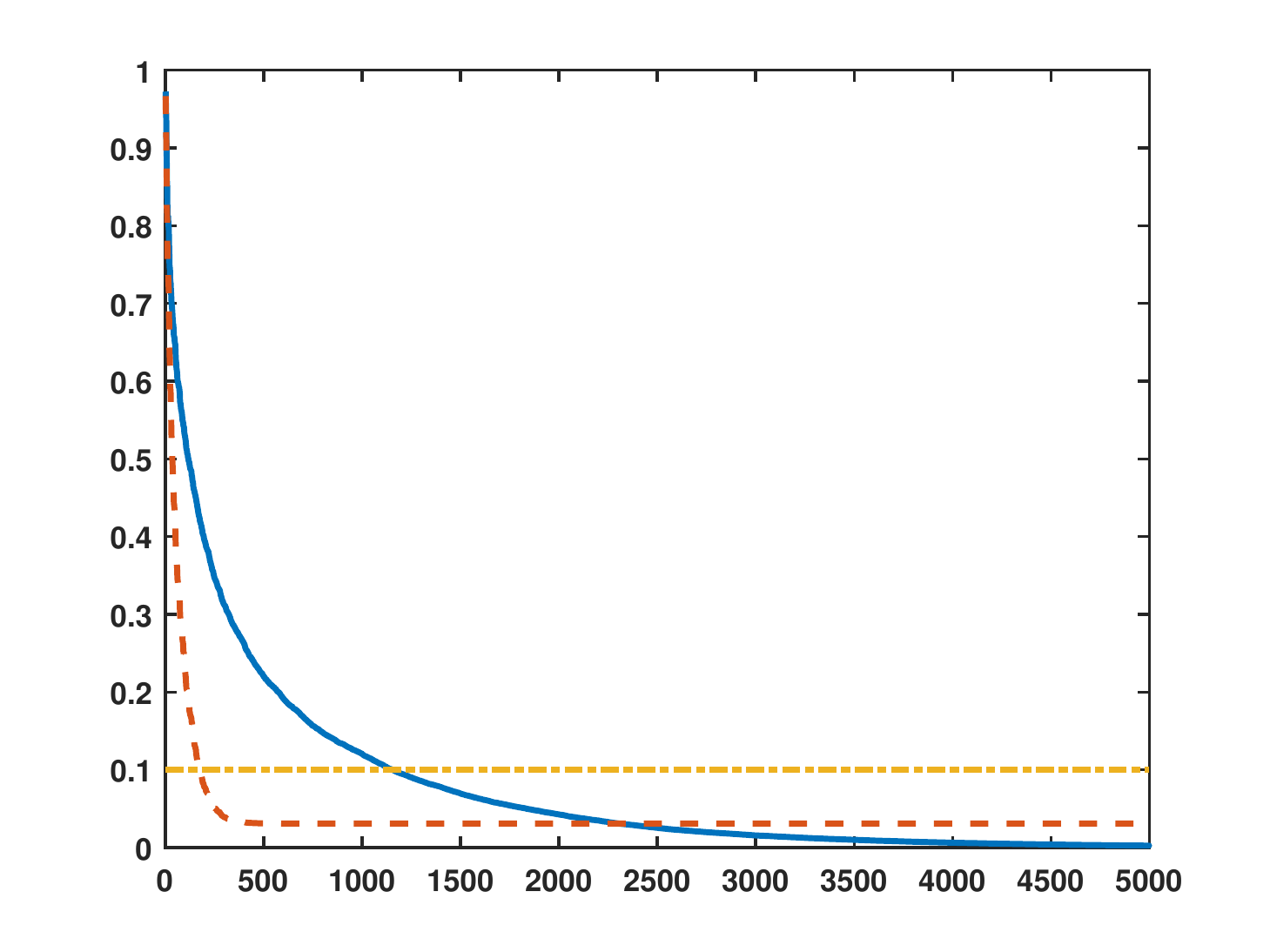}
\caption{The distance to the optimum~$\hat{x}$ for
the unregularized problem (shown in the solid blue line)
and the regularized problem (shown in the dashed orange line).
The desired regularization error bound~$\epsilon$ is shown
as the dash-dotted yellow line. We see that the regularized
problem converges substantially faster than the unregularized
one. Regularization error is equal to the final value of
the dashed orange line, which we see is indeed below the
dash-dotted yellow line, indicating that the desired
regularization error is indeed obeyed. 
}
\label{fig:network} 
\end{figure}

Both simulations were run for 2000 timesteps. 
To force asynchrony, 
at each timestep agent
$i$ has a 10\% chance of transmitting its state to agent $j$, and
agent $i$ has a 10\% chance of computing an update to its own state.

As shown in Figure~\ref{fig:network}, the regularized problem
(dashed orange line) shows notably faster
convergence than the unregularized one (solid blue line). 
The final error in the regularized solution was $e_{A}=0.0308$,
well below the desired upper bound of $e_{A}<0.1$, which
is shown in Figure~\ref{fig:network} by the dash-dotted yellow line. 
The condition number
of the regularized problem was $k_{Q+A}=8.1836$, satisfying our desired
condition number bound as well. 

\section{Conclusion} \label{sec:conclusion}
A totally asynchronous quadratic programming framework was presented.
This framework allowed agents to independently choose all stepsize
and regularization parameters and showed fast convergence
in simulation.  
Future work includes incorporating
functional constraints into this framework, as well as implementing this work
for path planning problems in robotic teams. 

\appendices{}

\section*{Appendix}

\textit{Proof of Theorem 1:} The following lemmas will facilitate the proof of
Theorem~\ref{thm:step}.

\begin{lemma} \label{lem:jordan}
Let $B$ and $C$ be positive definite $n\times n$
Hermitian matrices with eigenvalues $\lambda_{1}(B)\geq\lambda_{2}(B)\geq...\geq\lambda_{n}(B)$
and $\lambda_{1}(C)\geq\lambda_{2}(C)\geq...\geq\lambda_{n}(C)$.
Then 
\begin{equation}
\lambda_{n}(BC+CB)\geq\hspace{-0.05cm}\min_{\beta\in\{1,n\}}\hspace{-0.05cm}\left[\hspace{-0.05cm}\frac{\lambda_{\beta}(B)\lambda_{\beta}(C)}{2}\hspace{-0.05cm}\left(\hspace{-0.05cm}\frac{(\sqrt{k_{B}}\hspace{-0.05cm}+\hspace{-0.05cm}1)^{2}}{\sqrt{k_{B}}}\hspace{-0.05cm}-k_{C}\frac{(\sqrt{k_{B}}\hspace{-0.05cm}-\hspace{-0.05cm}1)^{2}}{\sqrt{k_{B}}}\hspace{-0.05cm}\right)\hspace{-0.05cm}\right],
\end{equation}

where $k_{B}$ and $k_{C}$ are the spectral condition numbers of
$B$ and $C$, respectively.
\end{lemma}

\textit{Proof:} See \cite{nicholson1979eigenvalue}. $\hfill\blacksquare$

\begin{lemma} \label{lem:boundswork}
Let $Q=Q^{T}\succ0$,~$Q \in \R^{n \times n}$ have
condition number $k_{Q}$, and let $\Gamma=\text{diag}\left(\gamma_{1}I_{n_{1}},...,\gamma_{N}I_{n_{N}}\right)$.
If
\begin{equation}
\gamma_{i}\in\left(\frac{\sqrt{k_{Q}}-1}{\|Q\|_{2}\sqrt{k_{Q}}},\frac{\sqrt{k_{Q}}+1}{\|Q\|_{2}\sqrt{k_{Q}}}\right) \text{ for all } i\in[N],
\end{equation}
then~$Q^{-1}\Gamma^{-1}+\Gamma^{-1}Q^{-1}-I\succ0$.
\end{lemma}

\textit{Proof:} Let $B=Q^{-1}$ and $C=\Gamma^{-1}$, and note that
\begin{equation}
k_{Q^{-1}}=\frac{\lambda_{1}(Q^{-1})}{\lambda_{n}(Q^{-1})}=\frac{\lambda_{n}^{-1}(Q)}{\lambda_{1}^{-1}(Q)}=\frac{\lambda_{1}(Q)}{\lambda_{n}(Q)}=k_{Q}
\end{equation}
and likewise $k_{\Gamma^{-1}}=k_{\Gamma}$. Now, using Lemma~\ref{lem:jordan}, we can write
\begin{equation} \label{eqn:bigjordan}
\lambda_{n}(Q^{-1}\Gamma^{-1}+\Gamma^{-1}Q^{-1})\geq\min_{\beta\in\{1,n\}}\Bigg[\frac{\lambda_{\beta}^{-1}(Q)\lambda_{\beta}^{-1}(\Gamma)}{2}\Bigg(\frac{(\sqrt{k_{Q}}+1)^{2}}{\sqrt{k_{Q}}} -k_{\Gamma}\frac{(\sqrt{k_{Q}}-1)^{2}}{\sqrt{k_{Q}}}\Bigg)\Bigg].
\end{equation}
Define the constants
\begin{equation}
\gamma_{lower}=\frac{\sqrt{k_{Q}}-1}{\|Q\|_{2}\sqrt{k_{Q}}} \text{  and  } \gamma_{upper}=\frac{\sqrt{k_{Q}}+1}{\|Q\|_{2}\sqrt{k_{Q}}}.
\end{equation}
Then, by hypothesis, $\gamma_{i}\in\left(\gamma_{lower},\gamma_{upper}\right)$ and
\begin{equation}
k_{\Gamma}<\frac{\gamma_{upper}}{\gamma_{lower}}=\frac{\sqrt{k_{Q}}+1}{\sqrt{k_{Q}}-1}.
\end{equation}
Substituting this bound into Equation~\eqref{eqn:bigjordan}. we find
\begin{align*}
 \lambda_{n}(Q^{-1}\Gamma^{-1}+\Gamma^{-1}Q^{-1}) & \geq\min_{\beta\in\{1,n\}}\Bigg[\frac{\lambda_{\beta}^{-1}(Q)\lambda_{\beta}^{-1}(\Gamma)}{2}\Bigg(\frac{(\sqrt{k_{Q}}+1)^{2}}{\sqrt{k_{Q}}}-\frac{(\sqrt{k_{Q}}+1)}{(\sqrt{k_{Q}}-1)}\frac{(\sqrt{k_{Q}}-1)^{2}}{\sqrt{k_{Q}}}\Bigg)\Bigg]\\
 & \geq\min_{\beta\in\{1,n\}}\left[\frac{\lambda_{\beta}^{-1}(Q)\lambda_{\beta}^{-1}(\Gamma)}{2}\left(2\frac{\sqrt{k_{Q}}+1}{\sqrt{k_{Q}}}\right)\right],
\end{align*}
the right hand side of which is always positive. This indicates that
the minimum will occur when $\beta=1$, and gives 
\begin{equation}
\lambda_{n}(Q^{-1}\Gamma^{-1}+\Gamma^{-1}Q^{-1})\geq\lambda_{1}^{-1}(Q)\lambda_{1}^{-1}(\Gamma)\left(\frac{\sqrt{k_{Q}}+1}{\sqrt{k_{Q}}}\right).
\end{equation}
Using $\lambda_{1}(Q)=\|Q\|_{2}$ and $\lambda_{1}(\Gamma)<\gamma_{upper}$,
we have
\begin{align*}
\lambda_{n}(Q^{-1}\Gamma^{-1}+\Gamma^{-1}Q^{-1}) & >\left(\frac{1}{\|Q\|_{2}}\right)\left(\frac{\|Q\|_{2}\sqrt{k_{Q}}}{\sqrt{k_{Q}}+1}\right)\left(\frac{\sqrt{k_{Q}}+1}{\sqrt{k_{Q}}}\right)\\
 & =1,
\end{align*}
which implies $Q^{-1}\Gamma^{-1}+\Gamma^{-1}Q^{-1}\succ I$. \hfill$\blacksquare$

From $Q^{-1}\Gamma^{-1}+\Gamma^{-1}Q^{-1}-I\succ0$
we see that the matrix on the left hand side is symmetric positive
definite. In particular, all eigenvalues are positive. Then, by Sylvester's
Law of Inertia \cite[Fact 5.8.17]{bernstein2009matrix}, performing a congruence transformation will 
preserve the positivity of the eigenvalues.
Any non-singular matrix $P$ provides a valid congruence transformation
of a matrix $M$ via $M'=P^{T}MP$, where $M'$ is the transformed
matrix. Using the nonsingular matrix $P=\Gamma Q$, we find
\begin{align*}
P^{T}(Q^{-1}\Gamma^{-1}+\Gamma^{-1}Q^{-1}-I)P & \succ0\\
Q\Gamma (Q^{-1}\Gamma^{-1}+\Gamma^{-1}Q^{-1}-I)\Gamma Q & \succ0,
\end{align*}
where we have used the symmetry of $Q$ and $\Gamma$. Expanding,
we find~$Q\Gamma+\Gamma Q-Q\Gamma^{2}Q\succ0$.
Multiplying by $-1$ and adding~$I$ gives
$I-Q\Gamma-\Gamma Q+Q\Gamma^{2}Q\prec I$,
which we factor as~$(I-\Gamma Q)^{T}(I-\Gamma Q)\prec I$.
This in turn implies~${\lambda_{max}[(I-\Gamma Q)^{T}(I-\Gamma Q)]<1}$,
where the square root gives
$\sqrt{\lambda_{max}[(I-\Gamma Q)^{T}(I-\Gamma Q)]}<1$
and finally~${\|I-\Gamma Q\|_{2}<1}$. \hfill $\blacksquare$

\bibliographystyle{IEEEtran}
\bibliography{Biblio}
\end{document}